\documentclass[a4paper,10pt]{article}
\usepackage[hmargin=2cm,vmargin=2cm]{geometry}
\usepackage{booktabs}
\usepackage{amsmath, amssymb, amsthm,xcolor}
\usepackage{enumerate}%

\usepackage{amsfonts}
\usepackage{float}
\usepackage[small,it]{caption}
\usepackage{multirow}
\usepackage{setspace}
\usepackage{pdfpages}
\usepackage{exscale,relsize}
\usepackage{subfig}
\usepackage{setspace}
\usepackage{bm}
\usepackage{empheq}
\usepackage{enumitem}
\usepackage[numbers,sort]{natbib}
\usepackage{comment}
\usepackage{graphicx}
\usepackage{epstopdf}
\newcount\colveccount
\newcommand*\colvec[1]{
        \global\colveccount#1
        \begin{pmatrix}
        \colvecnext
}
\def\colvecnext#1{
        #1
        \global\advance\colveccount-1
        \ifnum\colveccount>0
                \\
                \expandafter\colvecnext
        \else
                \end{pmatrix}
        \fi
}
\usepackage{ulem}
\numberwithin{equation}{section}
\theoremstyle{definition} 

\theoremstyle{plain} 

\theoremstyle{plain} 

\theoremstyle{definition} 

\theoremstyle{remark}

\theoremstyle{definition}

\begin{document}

\normalem
\title{Stability of travelling waves in a \emph{Wolbachia} invasion}
\author{Matthew H.\ Chan$^{1}$,
\and Peter S.\ Kim$^{2}$,
\and Robert Marangell$^{3}$}
\date{\today}
\maketitle

\footnotetext[1]{School of Mathematics and Statistics, University of Sydney,
 NSW 2006, Australia. {\tt M.Chan@maths.usyd.edu.au}}
   \footnotetext[2]{School of Mathematics and Statistics, University of Sydney,
 NSW 2006, Australia. {\tt pkim@maths.usyd.edu.au}}
  \footnotetext[3]{School of Mathematics and Statistics, University of Sydney,
 NSW 2006, Australia. {\tt robert.marangell@sydney.edu.au}}
\renewcommand{\thefootnote}{\arabic{footnote}}

\begin{abstract}
Numerous studies have examined the growth dynamics of \emph{Wolbachia} within populations and the resultant rate of spatial spread. This spread is typically characterised as a travelling wave with bistable local growth dynamics due to a strong Allee effect generated from cytoplasmic incompatibility. While this rate of spread has been calculated from numerical solutions of reaction-diffusion models, none have examined the spectral stability of such travelling wave solutions. In this study we analyse the stability of a travelling wave solution generated by the reaction-diffusion model of Chan \& Kim \cite{177} by computing the essential and point spectrum of the linearised operator arising in the model. The point spectrum is computed via an Evans function using the compound matrix method, whereby we find that it has no roots with positive real part. Moreover, the essential spectrum lies strictly in the left half plane. Thus, we find that the travelling wave solution found by Chan \& Kim \cite{177} corresponding to competition between \emph{Wolbachia}-infected and -uninfected mosquitoes is linearly stable. We employ a dimension counting argument to suggest that, under realistic conditions, the wavespeed corresponding to such a solution is unique.
\end{abstract}


\pagebreak

\section{Introduction}

\emph{Wolbachia} are common endosymbiotic bacteria that are estimated to infect up to 66\% of all insect species \cite{178}. They are primarily vertically transmitted and can induce reproductive phenotypes in its host to confer a reproductive advantage and hence improve chances of persistence. A well studied induced reproductive phenotype is \emph{cytoplasmic incompatibility} (CI), whereby offspring of infected males and uninfected females have an increased probability of death. This has lead to the proposal of a deliberate \emph{Wolbachia} introduction into wild mosquito populations to reduce transmission of vector-borne diseases, since particular CI-inducing \emph{Wolbachia} strains have been shown to reduce proliferation of various viruses (see Brelsfoard \& Dobson \cite{6} and references therein). Two particular CI-inducing strains of \emph{Wolbachia}, \emph{WMel} and \emph{WMelPop}, have received much attention due to evidence of these strains inhibiting dengue transmission in \emph{Aedes Aegypti} mosquitoes. 
\\
\\
Although CI has been proposed as a key mechanism for the success of \emph{Wolbachia}-based strategies, fitness reducing phenotypes are also a result of certain \emph{Wolbachia} strains. For example, both \emph{WMel} and \emph{WMelPop} reduce both fecundity and lifespan of infected females \cite{2,8,180}. Turelli \cite{3} first showed that the amalgamation of these \emph{Wolbachia}-induced effects result in a strong Allee effect for an invasion, that is, infection is only successful if initial infection densities are above the critical Allee threshold. This threshold has been found by various modelling studies to be dependent on the fecundity cost of infection, the strength of CI in reducing offspring from infected males and uninfected females, and the probability of successful \emph{Wolbachia} transmission \cite{179,4,7,13,177}. 
\\
\\
Many studies have analysed the growth dynamics of a deliberate \emph{Wolbachia} introduction via non-spatial models, but relatively few studies have studied the dynamics of the spatial spread (see Hancock \& Godfray \cite{181} and  Barton \& Turelli \cite{182} for examples). This study is motivated by the reaction-diffusion model of Chan \& Kim \cite{177}, who examine the spatial spread of \emph{Wolbachia} in a homogeneous environment by incorporating slow and fast compartments in their model. Chan \& Kim \cite{177} numerically show that there exists a travelling wave solution corresponding to a \emph{Wolbachia} invasion, that is, competition between \emph{Wolbachia}-infected and -uninfected mosquitoes, and estimate the wavespeed corresponding to this solution. A simplified model which uses a weighted average of the slow and fast diffusion coefficients was found to yield similar wavespeeds. This model assumes perfect vertical transmission of \emph{Wolbachia}, an increased fecundity for infected individuals, longevity reducing effects of $10\%$ and that \emph{Wolbachia} infection induces CI. Here we perform a spectral stability analysis by examining the linearised operator arising in this simplified model. We show that the essential spectrum is bounded to the left-half plane for all relevant biological parameter values and show that the point spectrum contains no elements in the right-half plane. Moreover, we show that there exists a travelling wave solution for only a unique wavespeed.

\section{Problem setup}
The non-dimensionalised system of partial differential equations from Chan \& Kim \cite{177} is given by

\begin{equation}
\label{PDEsystem}
\begin{cases}
u_t = u_{xx} - \rho u_x + u (1-S) - \alpha \mu u, \\
v_t = v_{xx} - \rho v_x + F v (1-S) (1-s_h A) - \mu v, \\
\end{cases}
\end{equation}
\\
where $u$ is the density of \emph{Wolbachia} infected mosquitoes, $v$ is the density of uninfected mosquitoes, $S=u+v$ and $A = u/S$. Following Chan \& Kim \cite{177}, we let $s_h=0.45$, $F = 1.0526$ and $\mu = 0.0162$, which correspond to parameter values for \emph{Aedes aegypti} at $30 ^\circ \text{C}$. The parameters $F^{-1}$, $\mu$, $\rho$ and $s_h$ correspond to the relative fecundity of infected females to uninfected females, mortality rate, advection rate and probability of embryo death due to cytoplasmic incompatibility, respectively. Additionally, we let $\alpha = 1.1$, which reflects the $10\%$ relative reduction in lifespan associated with the \emph{WMel} strain of \emph{Wolbachia}.
\\
\\
Converting to travelling wave coordinates $z = x-(c+\rho)t$ , we have that $(u(x,t),v(x,t)) = (u(z,t),v(z,t))$. This yields 

\begin{equation}
\label{travelsys}
\begin{cases}

u_t &= u_{zz} + c u_z + u (1-S)- \alpha \mu u,\\

v_t &= v_{zz} + c v_z + F v (1-S)(1-s_h A)-\mu v,\\

(u,v)(-\infty)&=\mathbf{e_{-}},\\

(u,v)(\infty)&=\mathbf{e_{+}},\\

(u',v')(\pm \infty) &= \mathbf{0}. 
\end{cases}
\end{equation}
\\
where $\mathbf{e_{-}}=(1-\alpha \mu,0)$ and $\mathbf{e_{+}}=(0,1-\frac{\mu}{F})$.
\\
\\
We linearise about the travelling wave solution $(\hat{u}(z),\hat{v}(z))$ via the substitution

\begin{equation}
\colvec{2}{u(z,t)}{v(z,t)} = \colvec{2}{\hat{u}(z)}{\hat{v}(z)} + \colvec{2}{p(z,t)}{q(z,t)},
\end{equation}
\\
where $p(z,t)$ and $q(z,t)$ are perturbations in $H^1$, $\forall t \in \mathbb{R}$. Collecting first order perturbation terms, we obtain

\begin{equation}
\label{travellingwavepde}
\begin{cases}

p_t &= p_{zz} + cp_z - \hat{u}(p+q) + p(1-\hat{S}) - \alpha \mu p,\\

q_t &= q_{zz} + cq_z + F \hat{v} \left( \frac{(\hat{S}^2-\hat{v})s_h}{\hat{S}^2} - 1 \right)p + F \left( 1-\hat{v} - \hat{S} + \frac{
\hat{u} s_h (\hat{S}^2-\hat{u})}{\hat{S}^2} \right)q - \mu q.\\
\end{cases}
\end{equation}
\\
We define the linear operator $\mathcal{L}$ by

\begin{equation}
\mathcal{L}\colvec{2}{p}{q} := \colvec{2}{p_{zz} + cp_z - \hat{u}(p+q) + p(1-\hat{S}) - \alpha \mu p}{q_{zz} + cq_z + F \hat{v} \left( \frac{(\hat{S}^2-\hat{v})s_h}{\hat{S}^2} - 1 \right)p + F \left( 1-\hat{v} - \hat{S} + \frac{\hat{u} s_h (\hat{S}^2-\hat{u})}{\hat{S}^2} \right)q - \mu q},
\end{equation}
\\
which has the corresponding eigenvalue problem $(\mathcal{L}-\lambda)\colvec{2}{p}{q}=0$. We introduce the substitutions $s=p_z$ and $t=q_z$ to convert the eigenvalue problem into a first order boundary value problem, and denote this equivalent operator of $(\mathcal{L}-\lambda)\colvec{2}{p}{q}=0$ by $\mathcal{T}(p,q,s,t)^T$, where $\mathcal{T}(\bm{y})(z) = \left( \frac{d}{dz} - A(z,\lambda) \right) \bm{y}$ and $\bm{y} = (p,q,s,t)^T$. This process yields

\begin{equation}
\label{bvp}
\begin{cases}
\begin{aligned}
\bm{y}'(z) &= A(z,\lambda) \bm{y}(z),\\
\bm{y}(-\infty) &= \left( 1- \alpha \mu,0,0,0 \right)^T,\\
\bm{y}(\infty) &= \left( 0,1-\frac{\mu}{F},0,0 \right)^T,
\end{aligned}
\end{cases}
\end{equation}
\\
and

\begin{equation}
\label{nonconmat}
A(z,\lambda) := \left( \begin{array}{cccc} 0 & 0 & 1 & 0 \\ 0 & 0 & 0 & 1 \\ \hat{u}-(1-\hat{S})+\alpha \mu+\lambda & \hat{u} & -c & 0 \\ F \hat{v}\left(1- \frac{(\hat{S}^2-\hat{v})s_h}{\hat{S}^2} \right) & F \left( -1+\hat{v} + \hat{S} - \frac{u s_h (\hat{S}^2-\hat{u})}{\hat{S}^2} \right) +\mu +\lambda & 0 & -c \end{array} \right).
\end{equation}
\\
To assess the stability of travelling wave solution $\hat{u}$ and $\hat{v}$ we need to locate the spectrum of the linearised operator $\mathcal{L}$ as an operator on $H^1 \times H^1$. If $(\mathcal{L} - \lambda)^{-1}$ does not exist or is unbounded for $\lambda \in \mathbb{C}$, then $\lambda$ is in the \emph{spectrum} $\sigma(\mathcal{L})$ of the operator $\mathcal{L}$. The complement of the spectrum in $\mathbb{C}$ is the \emph{resolvent set} of $\mathcal{L}$. Following Kapitula \& Promislow \cite[Section 2.2.5]{170}, we define $\text{ind}(\mathcal{L}) = \dim[\ker (\mathcal{L})]-\text{codim}[\text{R}(\mathcal{L})]$ as the Fredholm index of $\mathcal{L}$, where $\text{R}(\mathcal{L})$ denotes the range of $\mathcal{L}$. The spectrum of a Fredholm operator $\mathcal{L}$ is decomposed into two sets:

\begin{enumerate}[label=(\roman*)]
  \item The essential spectrum, defined by
  $$
  \sigma_{\textnormal{ess}}(\mathcal{L}) = \{  \lambda \in \mathbb{C} \mid \lambda - \mathcal{L} \text{ is not Fredholm or } \lambda  - \mathcal{L} \text{ is Fredholm, but } \text{ind}(\lambda - \mathcal{L}) \neq 0 \}.
  $$
  \item The point spectrum, defined by
  $$
  \sigma_{\textnormal{pt}}(\mathcal{L}) = \{  \lambda \in \mathbb{C} \mid \text{ind} (\lambda - \mathcal{L}) = 0 \text{, but } \lambda - \mathcal{L} \text{ is not invertible}  \}.
  $$
\end{enumerate}

\noindent We define $A_\pm (\lambda) := \lim_{z\to \pm \infty} A(z,\lambda)$, which are given by

\begin{equation}
\label{Aminus}
A_{-}(\lambda) = \left( \begin{array}{cccc} 0 & 0 & 1 & 0 \\ 0 & 0 & 0 & 1 \\ 1- \alpha \mu+\lambda & 1- \alpha \mu & -c & 0 \\ 0 & -F \alpha \mu (1-s_h)+\mu+\lambda & 0 & -c \end{array} \right),
\end{equation}

\begin{equation}
\label{Aplus}
A_{+}(\lambda) = \left( \begin{array}{cccc} 0 & 0 & 1 & 0 \\ 0 & 0 & 0 & 1 \\ \mu \left( \alpha-\frac{1}{F} \right)+\lambda & 0 & -c & 0 \\ F-\mu+\mu s_h & F-\mu+\lambda & 0 & -c \end{array} \right).
\end{equation}
\\
The asymptotic operator of $\mathcal{T}$ is given by

\begin{equation}
\mathcal{T_\infty}(\lambda)\colvec{4}{p}{q}{s}{t} := \colvec{4}{p}{q}{s}{t}' - A_\infty (\lambda) \colvec{4}{p}{q}{s}{t},
\end{equation}
\\
where $A_\infty$ is the piecewise spatially constant matrix

\begin{equation}
A_\infty(\lambda) =
\begin{cases}
A_-(\lambda) \quad\quad z < 0, \\
A_+(\lambda) \quad\quad z \geq 0.
\end{cases}
\end{equation}

\section{Essential spectrum}

\begin{figure}[H]
  \centerline{\includegraphics[scale=0.5]{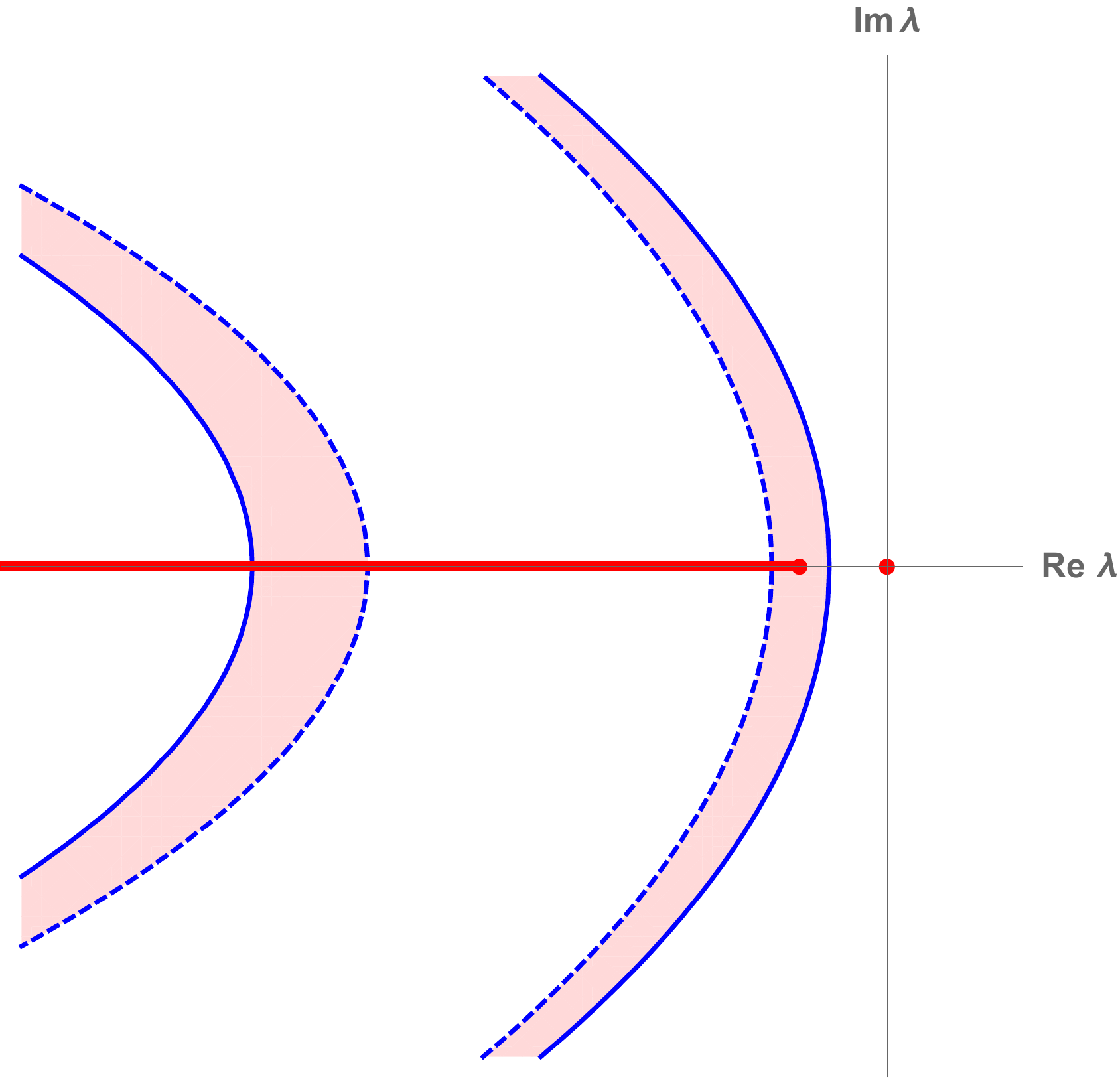}}
  \caption{The essential spectrum of $\mathcal{L}$ is given by $\lambda$ in the shaded regions. The blue dashed and solid lines represent $\lambda_-^{1,2}$ and $\lambda_+^{1,2}$ respectively. The red line indicates the absolute spectrum given by Eq.~(\ref{absspectrum}) and the red dot at the origin represents an eigenvalue. Note that these are not drawn to scale for visualisation purposes.}
  \label{essentialspectrumplot}
\end{figure}

The PDE given in (\ref{travelsys}) is autonomous, and so the only non-constant coefficients are due to the functions $\hat{u}(z)$ and $\hat{v}(z)$ in the reaction terms. These are heteroclinic orbits in phase space (connecting  $\mathbf{e_{-}}=(1-\alpha \mu,0)$ and $\mathbf{e_{+}}=(0,1-\frac{\mu}{F})$) which decay exponentially as $z \to \pm \infty$ as shown in Section \ref{waveprofilesect}. This shows that $\mathcal{L}$ is exponentially asymptotic. From Kapitula \& Promislow \cite[Theorem 3.1.11]{170} it follows that $\mathcal{L}$ is a relatively compact perturbation of the asymptotic operator $\mathcal{L}_\infty$, {defined as the limit of $\mathcal{L}$ as $z \to \pm \infty$ and equivalent to the operator $\mathcal{T}_\infty(0)$}. Then by Weyl's Essential Spectrum Theorem, we have that $\sigma_{\textnormal{ess}}(\mathcal{L}) = \sigma_{\textnormal{ess}}(\mathcal{L}_\infty)$, or equivalently $\sigma_{\textnormal{ess}}(\mathcal{T}) = \sigma_{\textnormal{ess}}(\mathcal{T}_\infty)$ \cite[Theorem 2.2.6]{170}.
\\
\\
A crucial concept behind the spectrum of an operator is the existence of an exponential dichotomy. Essentially, this states that each solution to (\ref{bvp}) decays exponentially either in forward or backward $z$. For spatially constant matrices, the existence of an exponential dichotomy simply means that the matrix is hyperbolic. We define the \emph{Morse index} of a constant matrix $A$ to be the dimension of the unstable subspace associated with $A$, and let $i_\pm(\lambda)$ denote the Morse indices of the asymptotic matrices $A_\pm$, given by Eq.~(\ref{Aminus}) {and} Eq.~({\ref{Aplus}}). It can be shown that for $\lambda \in \mathbb{C}$ such that $\mathcal{T}_\infty$ is Fredholm, we have $\textnormal{ind}(\mathcal{T}_\infty-\lambda)={i_-(\lambda) - i_+(\lambda)}$ \cite[Lemma 3.1.10]{170}. Thus, we can characterise the essential spectrum of $\mathcal{L_\infty}$ as

\begin{equation}
\sigma_{\textnormal{ess}}(\mathcal{L}_\infty) = \left\{ \lambda \in \mathbb{C} \mid i_-(\lambda) \neq i_+(\lambda) \right\} \cup \left\{ \lambda \in \mathbb{C} \mid \textnormal{dim}~ \mathbb{E}^\textnormal{c}(A_\pm(\lambda)) \neq 0 \right\},
\end{equation}

\noindent where $\mathbb{E}^c$ denotes the center {subspace} associated with the {asymptotic} linearised system. 
\\
\\ 
The spatial eigenvalues of $A_-(\lambda)$ and $A_+(\lambda)$ are respectively given by

\begin{equation}
\label{eig1}
\eta_{-} = \frac{1}{2} \left( -c \pm \sqrt{c^2 + 4(1+\lambda- \alpha \mu)} \right), \frac{1}{2} \left( -c \pm \sqrt{c^2 + 4(\lambda+\mu(1- \alpha F(1-s_h))} \right)
\end{equation}
\\
and 

\begin{equation}
\label{eig2}
\eta_{+} = \frac{1}{2} \left( -c \pm \sqrt{c^2 + 4(F+\lambda-\mu)} \right), \frac{1}{2} \left( -c \pm \sqrt{c^2 + 4 \left( \lambda+\mu \left(\alpha-\frac{1}{F} \right) \right)} \right).
\end{equation}
\\
These spatial eigenvalues are non-hyperbolic when $\eta_\pm = ik$. Substituting this into Eq.~(\ref{eig1}) and (\ref{eig2}) respectively yields the {\em dispersion relations} 

\begin{equation}
\label{set1lambda}
\lambda_-^{1,2}(k) = -1+ \alpha \mu - k^2 + ick, -\mu (1- \alpha F(1-s_h)) - k^2 + ick,
\end{equation}

\begin{equation}
\label{set2lambda}
\lambda_+^{1,2}(k) = -F+\mu-k^2 + ick, \mu \left( \frac{1}{F} - \alpha \right) -k^2 + ick.
\end{equation}
\\
These form four parabolas in the complex plane parametrised by $k$. For $\lambda$ in between the region bounded by $\lambda_-^{1,2}(k)$, $A_-(\lambda)$ has three stable eigenvalues and one unstable eigenvalue; to the left of the region $A_-(\lambda)$ has four stable eigenvalues and to the right of the region $A_-(\lambda)$ has two stable and two unstable eigenvalues. This is also true for $\lambda_+^{1,2}(k)$ and $A_+(\lambda)$. Thus the essential spectrum is given by the region bounded between $\lambda_-^1$ and $\lambda_+^1$, and also between $\lambda_-^2$ and $\lambda_+^2$; this is shown in Figure \ref{essentialspectrumplot}.

\noindent 
{It will be convenient for us later on to know the location of the so-called {\em absolute spectrum}. The absolute spectrum is not spectrum per se, but its location characterises the breakdown of the analytic continuation (in terms of the spectral parameter $\lambda$) of the stable and unstable eigenspaces of the matrices $A_\pm(\lambda)$. It thus follows that the absolute spectrum coincides with a branch cut of the Evans function \cite[Section 3.2]{170}. For the case at hand, the absolute spectrum can be defined \cite{SS00} as the set in the complex plane where a pair of the eigenvalues of $A_+(\lambda)$ have equal real parts. This is the set


\begin{equation}
\label{absspectrum}
\left \{ \lambda \in \mathbb{R} \, \middle | \, \lambda \leq \mu \left( \frac{1}{F} - \alpha \right) - \frac{c^2}{4} \right  \}.
\end{equation}}
\\
We note that from Eq.~(\ref{absspectrum}), (\ref{set1lambda}) and (\ref{set2lambda}), the continuous and absolute spectrum are always bounded to the left-half plane for biologically relevant parameter constraints $F>1$, $\mu>0$, $s_h>0$ and $\alpha>1$.

\section{Point spectrum}

The existence of a travelling wave solution implies a heteroclinic connection in (\ref{travellingwavepde}) between the equilibria $\mathbf{e_{-}}=(1-\alpha \mu,0)$ and $\mathbf{e_{+}}=(0,1-\frac{\mu}{F})$. We denote the unstable subspace of the matrix $A_-$ by $U_-$ and the stable subspace of the matrix $A_+$ by $S_+$. To the right of the essential spectrum, we have that the dimension of $U_-$, which we denote by $k$, and dimension of $S_+$ sum to $4$, the dimension of the entire phase space. For our case, $k=2$.
\\
\\
The unstable eigenvalues of $A_-$ are given by

\begin{equation}
\label{unstableeigUm}
\eta^-_{1,2} = \frac{1}{2} \left( -c + \sqrt{c^2 + 4(1+\lambda-\alpha\mu)} \right), \frac{1}{2} \left( -c + \sqrt{c^2 + 4(\lambda+\mu(1- \alpha F(1-s_h))} \right)
\end{equation}
\\
and the stable eigenvalues of $A_+$ are given by

\begin{equation}
\label{unstableeigUp}
\eta^+_{1,2} = \frac{1}{2} \left( -c - \sqrt{c^2 + 4(F+\lambda-\mu)} \right), \frac{1}{2} \left( -c - \sqrt{c^2 + 4 \left( \lambda+\mu \left(\alpha-\frac{1}{F} \right) \right)} \right),
\end{equation}
\\
where we have $\eta^-_1 > \eta^-_2 > 0$ and $\eta^+_1 < \eta^+_2 < 0$. We denote $\boldsymbol{\zeta}^-_{1,2}$ and $\boldsymbol{\zeta}^+_{1,2}$ as the eigenvectors corresponding to $\eta^-_{1,2}$ and $\eta^+_{1,2}$ respectively.
\\
\\
We initialise Eq.~(\ref{bvp}) at $z=-\infty$ with $\boldsymbol{\zeta}^-_{1}$, $\boldsymbol{\zeta}^-_{2}$ and at $z=\infty$ with $\boldsymbol{\zeta}^+_{1}$, $\boldsymbol{\zeta}^+_{2}$, and solve the system towards a matching point, which we pick to be $z=0$. We denote the solutions of the former by $\mathbf{w}^-_1(0,\lambda)$, $\mathbf{w}^-_2(0,\lambda)$ and the latter by $\mathbf{w}^+_1(0,\lambda)$, $\mathbf{w}^+_2(0,\lambda)$, where $\mathbf{w}^-_i(z,\lambda)$ satisfy

\begin{equation}
\begin{cases}
\begin{aligned}
\frac{d}{dz} \mathbf{w}^-_i(z,\lambda) &= A(z,\lambda) \mathbf{w}^-_i(z,\lambda),\\
\mathbf{w}^-_i(z,\lambda) &\sim \exp(\eta^-_i z)\boldsymbol{\zeta}^-_i \quad\quad \text{for } \quad z \ll 0
\end{aligned}
\end{cases}
\end{equation}
\\
and $\mathbf{w}^+_i(z,\lambda)$ satisfy

\begin{equation}
\begin{cases}
\begin{aligned}
\frac{d}{dz} \mathbf{w}^+_i(z,\lambda) &= A(z,\lambda) \mathbf{w}^+_i(z,\lambda),\\
\mathbf{w}^+_i(z,\lambda) &\sim \exp(\eta^+_i z)\boldsymbol{\zeta}^+_i \quad\quad \text{for } \quad z \gg 0.
\end{aligned}
\end{cases}
\end{equation}
\\
The Evans function is defined by

\begin{equation}
D(\lambda) = \det \left[ \mathbf{w}^-_1(0,\lambda) , \mathbf{w}^-_2(0,\lambda), \mathbf{w}^+_1(0,\lambda) , \mathbf{w}^+_2(0,\lambda) \right],
\end{equation}
\\
which has the property that $D(\lambda)=0$ if and only if $\lambda$ is in the point spectrum of the operator $\mathcal{L}$. Roots of the Evans function correspond to solutions of the boundary value problem defined by Eq.~(\ref{bvp}), which decay appropriately as $z \to \pm \infty$.

\subsection{Compound matrix method}
The Evans function is numerically difficult to compute due to the stiffness of the problem stemming from the difficulty of resolving different modes of growth and decay. For example, since we have that $\eta^-_1 > \eta^-_2$, any numerical errors occurring when solving for $\mathbf{w}^-_2(0,\lambda)$ will grow at a rate proportionate to $\exp(\eta^-_1 z)$. Thus, although the solutions $\mathbf{w}^-_1(z,\lambda)$, $\mathbf{w}^-_2(z,\lambda)$ are linearly independent at $z=-\infty$, they quickly become numerically linearly dependent. Several methods have been proposed to overcome this issue such as the method of continuous orthogonalisation, the compound matrix method, Magnus methods and Grassmanian spectral shooting \cite{173,174,175,176,185}. Following Allen \& Bridges \cite{172}, we employ the compound matrix method which converts the problem into the six-dimensional wedge product space $\wedge^2 (\mathbb{C}^4)$, with basis $B = \{ \mathbf{e}_1 \wedge \mathbf{e}_2 , \mathbf{e}_1 \wedge \mathbf{e}_3 , \mathbf{e}_1 \wedge \mathbf{e}_4 , \mathbf{e}_2 \wedge \mathbf{e}_3 , \mathbf{e}_2 \wedge \mathbf{e}_4 , \mathbf{e}_3 \wedge \mathbf{e}_4\}$, where $\{ \mathbf{e}_1 , \mathbf{e}_2 , \mathbf{e}_3 , \mathbf{e}_4 \}$ is the standard basis for $\mathbb{C}^4$. The numerical advantage of this approach is that the evolution of $\mathbf{w}^-_{1}$, $\mathbf{w}^-_{2}$ and $\mathbf{w}^+_{1}$, $\mathbf{w}^+_{2}$ are incorporated into a single trajectory given by $\mathbf{w}^-_1 \wedge \mathbf{w}^-_2$ and $\mathbf{w}^+_1 \wedge \mathbf{w}^+_2$ respectively. The coordinate vector of $\mathbf{w}^-_1 \wedge \mathbf{w}^-_2$ relative to the basis $B$ is given by $[\mathbf{w}^-_1 \wedge \mathbf{w}^-_2]_B := \boldsymbol{\phi^-} = (\phi^-_1,\phi^-_2,\phi^-_3,\phi^-_4,\phi^-_5,\phi^-_6)$, where

\begin{equation}
\label{minordefn}
\begin{aligned}
\phi^-_{1} &=
\begin {vmatrix}
\mathbf{w}^-_{1,1} & \mathbf{w}^-_{2,1}\\
\mathbf{w}^-_{1,2} & \mathbf{w}^-_{2,2}
\end{vmatrix},
\hspace{0.2in}
\phi^-_{2} =
\begin {vmatrix}
\mathbf{w}^-_{1,1} & \mathbf{w}^-_{2,1}\\
\mathbf{w}^-_{1,3} & \mathbf{w}^-_{2,3}
\end{vmatrix},
\hspace{0.2in}
\phi^-_{3} =
\begin {vmatrix}
\mathbf{w}^-_{1,1} & \mathbf{w}^-_{2,1}\\
\mathbf{w}^-_{1,4} & \mathbf{w}^-_{2,4}
\end{vmatrix},
\\
\\
\phi^-_{4} &=
\begin {vmatrix}
\mathbf{w}^-_{1,2} & \mathbf{w}^-_{2,2}\\
\mathbf{w}^-_{1,3} & \mathbf{w}^-_{2,3}
\end{vmatrix},
\hspace{0.2in}
\phi^-_{5} =
\begin {vmatrix}
\mathbf{w}^-_{1,2} & \mathbf{w}^-_{2,2}\\
\mathbf{w}^-_{1,4} & \mathbf{w}^-_{2,4}
\end{vmatrix},
\hspace{0.2in}
\phi^-_{6} =
\begin {vmatrix}
\mathbf{w}^-_{1,3} & \mathbf{w}^-_{2,3}\\
\mathbf{w}^-_{1,4} & \mathbf{w}^-_{2,4}
\end{vmatrix}
\end{aligned}
\end{equation}
\\
and the second subscript in $\mathbf{w}^-_{i,j}$ denotes the $j$th element within the vector $\mathbf{w}^-_{i}$. Similarly, $[\mathbf{w}^+_1 \wedge \mathbf{w}^+_2]_B := \boldsymbol{\phi}^+$ is given by Eq.~(\ref{minordefn}) with $\mathbf{w}^-_{i,j}$ replaced by $\mathbf{w}^+_{i,j}$.
\\
\\
It can be shown (see Allen \& Bridges \cite{172}) that $\boldsymbol{\phi}(z) = \boldsymbol{\phi^-}(z),\boldsymbol{\phi^+}(z)$ satisfy the equation

\begin{equation}
\label{compoundmatIVP}
\boldsymbol{\phi}' = \tilde{A}(z,\lambda) \boldsymbol{\phi},
\end{equation}
\\
where $\tilde{A}$ is the induced matrix given by

\begin{equation}
\tilde{A} =
\begin{pmatrix}
A_{11}+A_{22} & A_{23} & A_{24} & -A_{13} & -A_{14} & 0 \\ 
A_{32} & A_{11}+A_{33} & A_{34} & A_{12} & 0 & -A_{14} \\ 
A_{42} & A_{43} & A_{11} + A_{44} & 0 & A_{12} & A_{13} \\ 
-A_{31} & A_{21} & 0 & A_{22} + A_{33} & A_{34} & -A_{24} \\ 
-A_{41} & 0 & A_{21} & A_{43} & A_{22}+A_{44} & A_{23} \\ 
0 & -A_{41} & A_{31} & -A_{42} & A_{32} & A_{33} + A_{44} 
\end{pmatrix},
\end{equation}
\\
with $A_{ij}$ given by Eq.~(\ref{nonconmat}).
\\
\\
We initialise the problem at $z=-\infty$ with $\boldsymbol{\phi}^-(-\infty)$ and at $z=\infty$ with $\boldsymbol{\phi}^+(\infty)$ and solve towards the matching point $z=0$. The Evans function is then defined to be

\begin{equation}
\begin{aligned}
\label{evansfunccompound}
D(\lambda) &= \mathbf{w}^-_1 \wedge \mathbf{w}^-_2 \wedge \mathbf{w}^+_1 \wedge \mathbf{w}^+_2, \\
&= \phi^-_1 \phi^+_6 - \phi^-_2 \phi^+_5 + \phi^-_3 \phi^+_4 + \phi^-_4 \phi^+_3 - \phi^-_5 \phi^+_2 + \phi^-_6 \phi^+_1.
\end{aligned}
\end{equation}
\\
For numerical stability, we scale the solution according to its exponential growth and decay rates by letting

\begin{equation}
\begin{cases}
\begin{aligned}
\boldsymbol{\phi}^-(z) &= \boldsymbol{\psi}^-(z) e^{-(\eta^-_1 + \eta^-_2)z},\\
\boldsymbol{\phi}^+(z) &= \boldsymbol{\psi}^+(z) e^{-(\eta^+_1 + \eta^+_2)z},
\end{aligned}
\end{cases}
\end{equation}
\\
which leads to

\begin{equation}
\begin{cases}
\begin{aligned}
\frac{d \boldsymbol{\psi}^-(z)}{dz} &= (\tilde{A} - (\eta^-_1 + \eta^-_2) \mathbb{I})\boldsymbol{\psi}^-(z) \quad\quad \text{for} \quad\quad z<0,\\
\frac{d \boldsymbol{\psi}^+(z)}{dz} &= (\tilde{A} - (\eta^+_1 + \eta^+_2) \mathbb{I})\boldsymbol{\psi}^+(z) \quad\quad \text{for} \quad\quad z>0,
\end{aligned}
\end{cases}
\end{equation}
\\
where the Evans function is equivalent to Eq.~(\ref{evansfunccompound}) except with $\phi$ replaced with $\psi$.
\\
\\
Since the Evans function is analytic to the right of the essential spectrum, we have via the Argument Principle that

\begin{equation}
\frac{1}{2\pi i}\oint_C \frac{D'(\lambda)}{D(\lambda)} \, d \lambda = N,
\end{equation}
\\
where $N$ is the number of zeroes in the interior of the region enclosed by $C$.
\\
\\
To check for eigenvalues of $\mathcal{L}$ with positive real part, we set up a closed semi-circle contour $C$ excluding the origin, as shown below in Figure \ref{contourplot}. We let $r_s$ and $r_b$ denote the radius of the smaller and larger circular arc respectively.

\begin{figure}[H]
  \centerline{\includegraphics[scale=0.5]{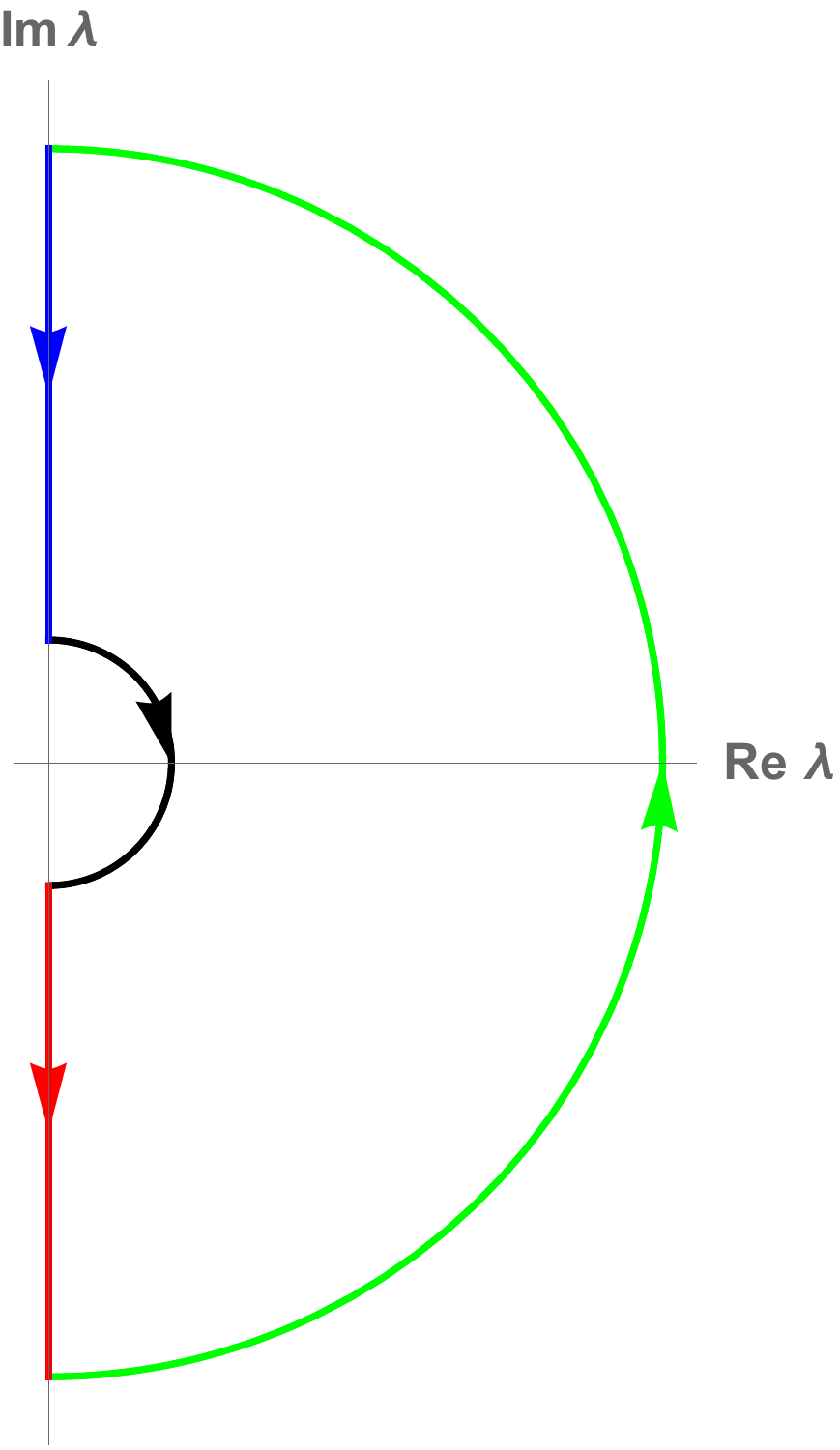}}
  \caption{The contour $C$.}
  \label{contourplot}
\end{figure}

\begin{figure}[H]
  \centerline{\includegraphics[scale=0.3]{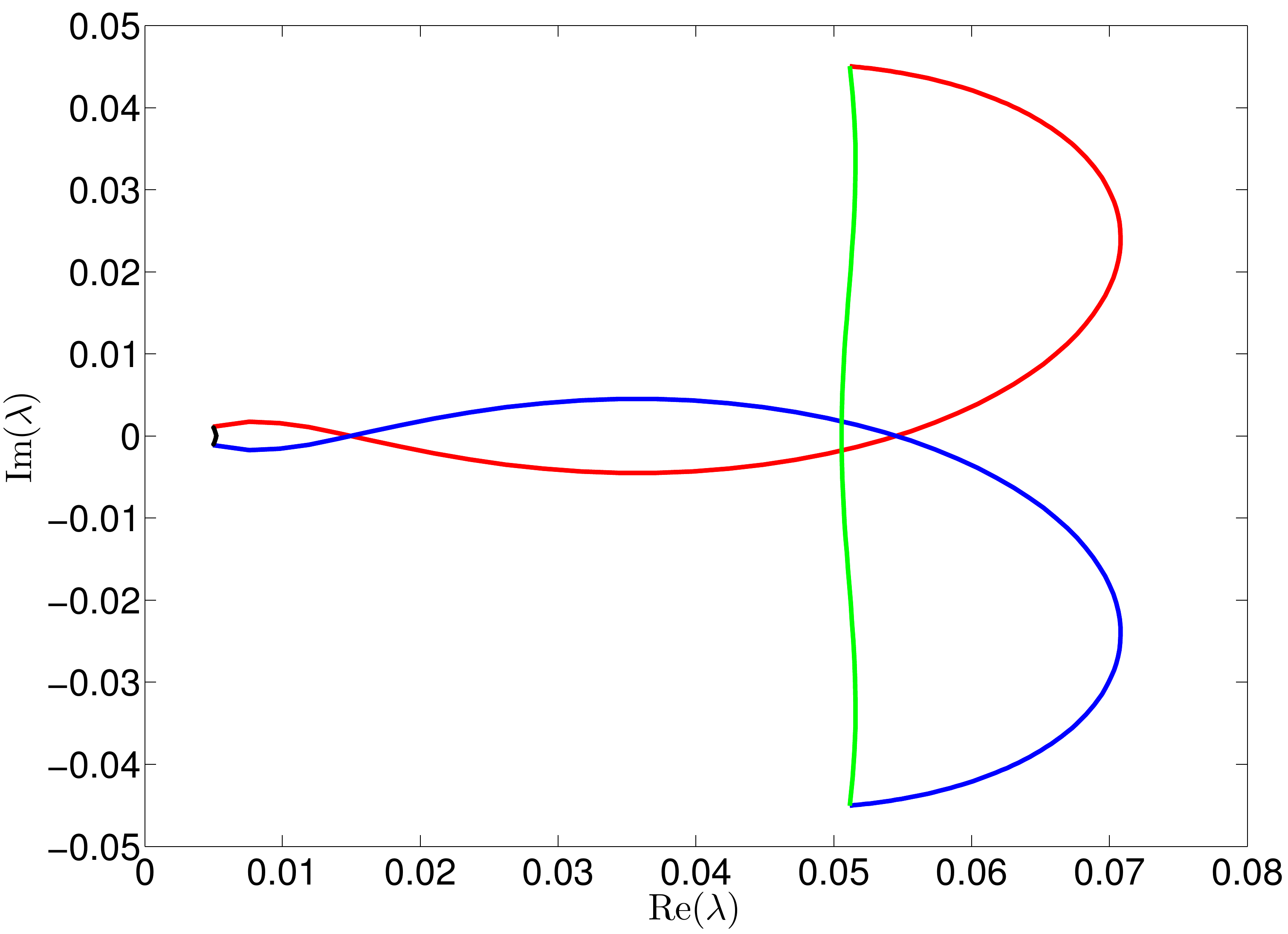}}
  \caption{The image of $C$ under $D(\lambda)$, where $r_s = 0.1$ and $r_b = 10$.}
  \label{evanscontoursmall}
\end{figure}

\begin{figure}[H]
  \hspace{1mm}\centerline{
  \subfloat[]{\label{evanscontourbig1}\includegraphics[scale=0.3]{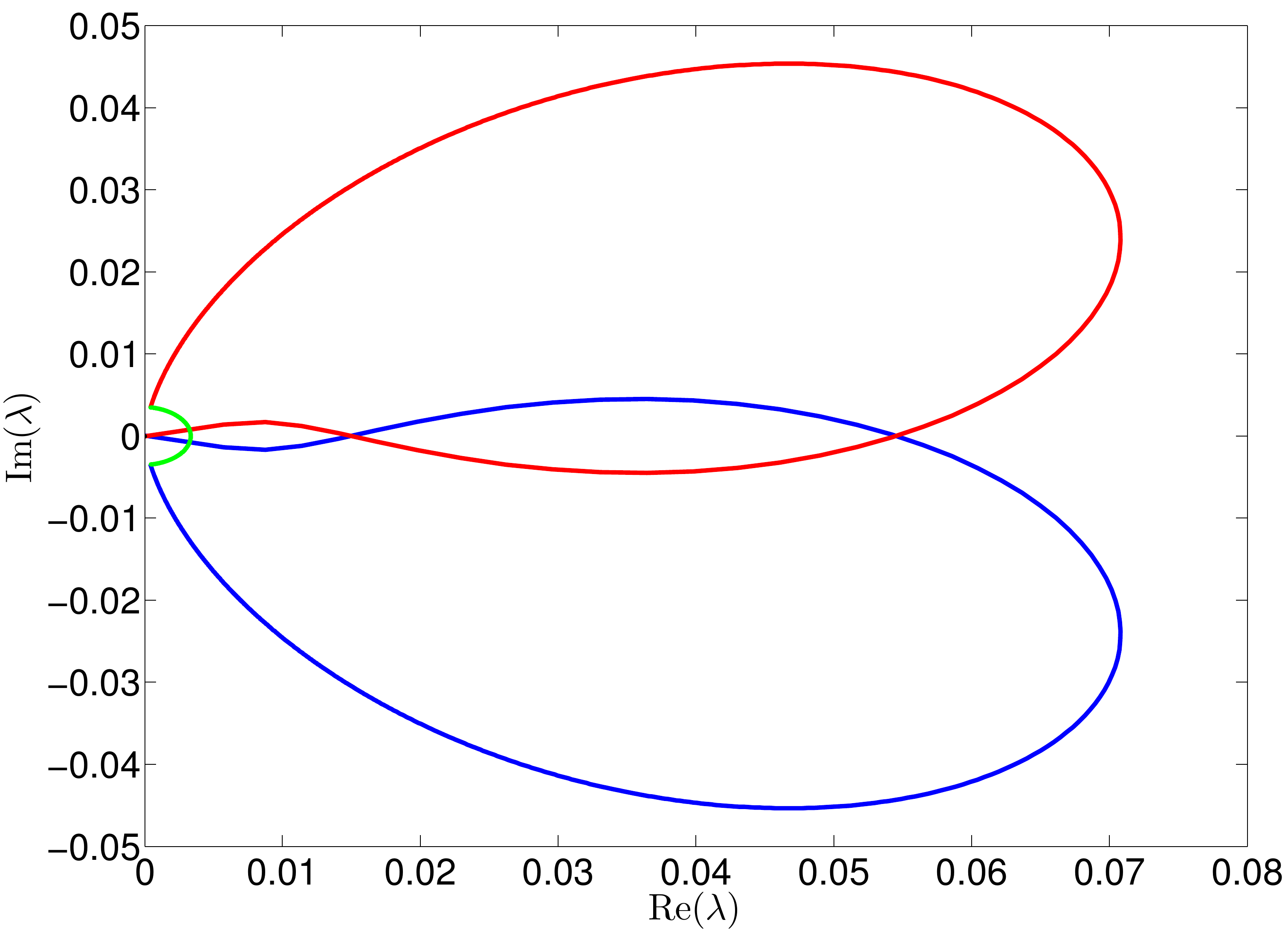}}
  \hspace{1mm}
  \subfloat[]{\label{evanscontourbig2}\includegraphics[scale=0.3]{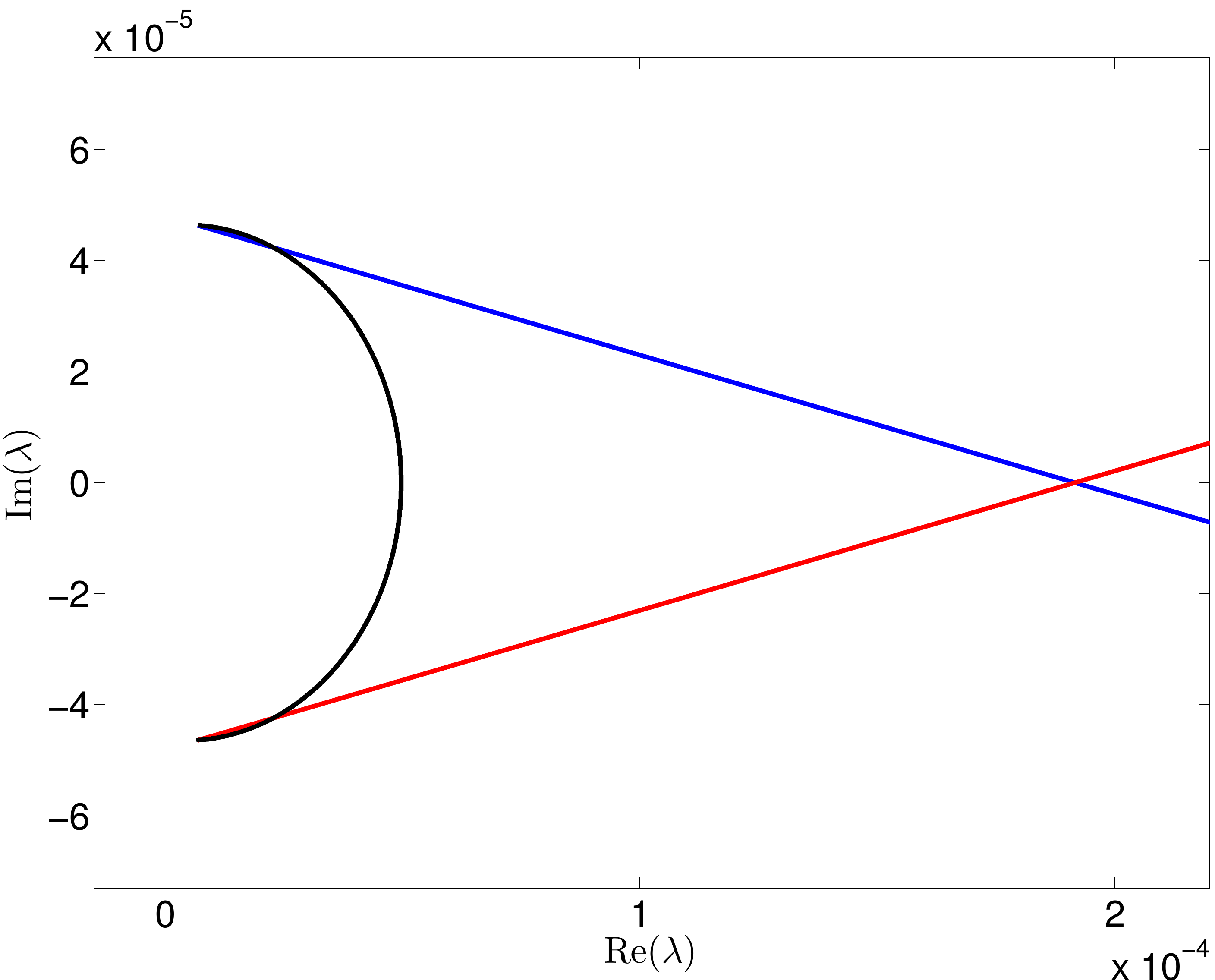}}}
  \caption{The image of $C$ under $D(\lambda)$, where $r_s = 0.001$ and $r_b = 500$.}
  \label{evanscontourbig}
\end{figure}

\begin{figure}[H]
  \hspace{1mm}\centerline{
  \subfloat[]{\label{argsmall}\includegraphics[scale=0.3]{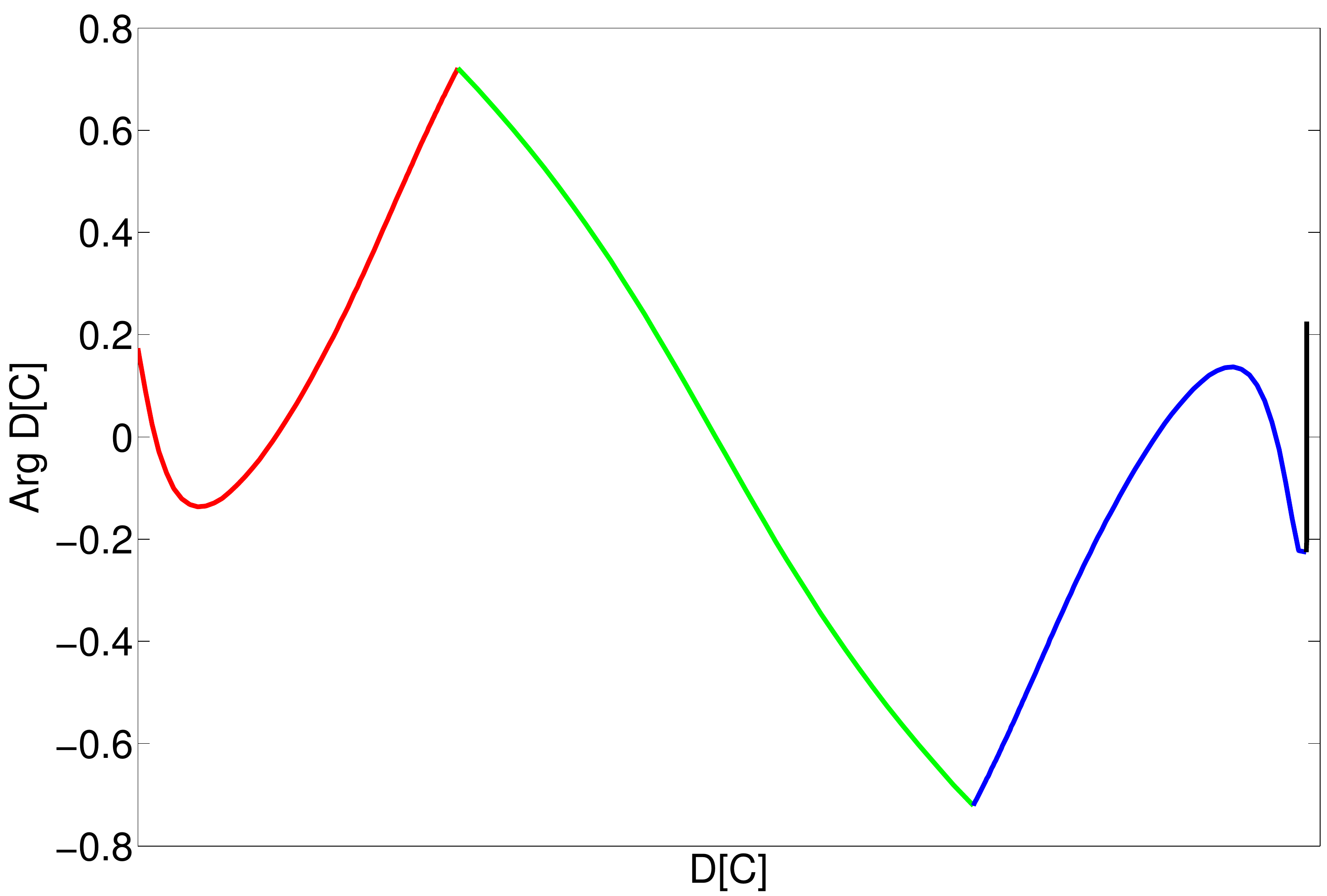}}
  \hspace{1mm}
  \subfloat[]{\label{argbig}\includegraphics[scale=0.3]{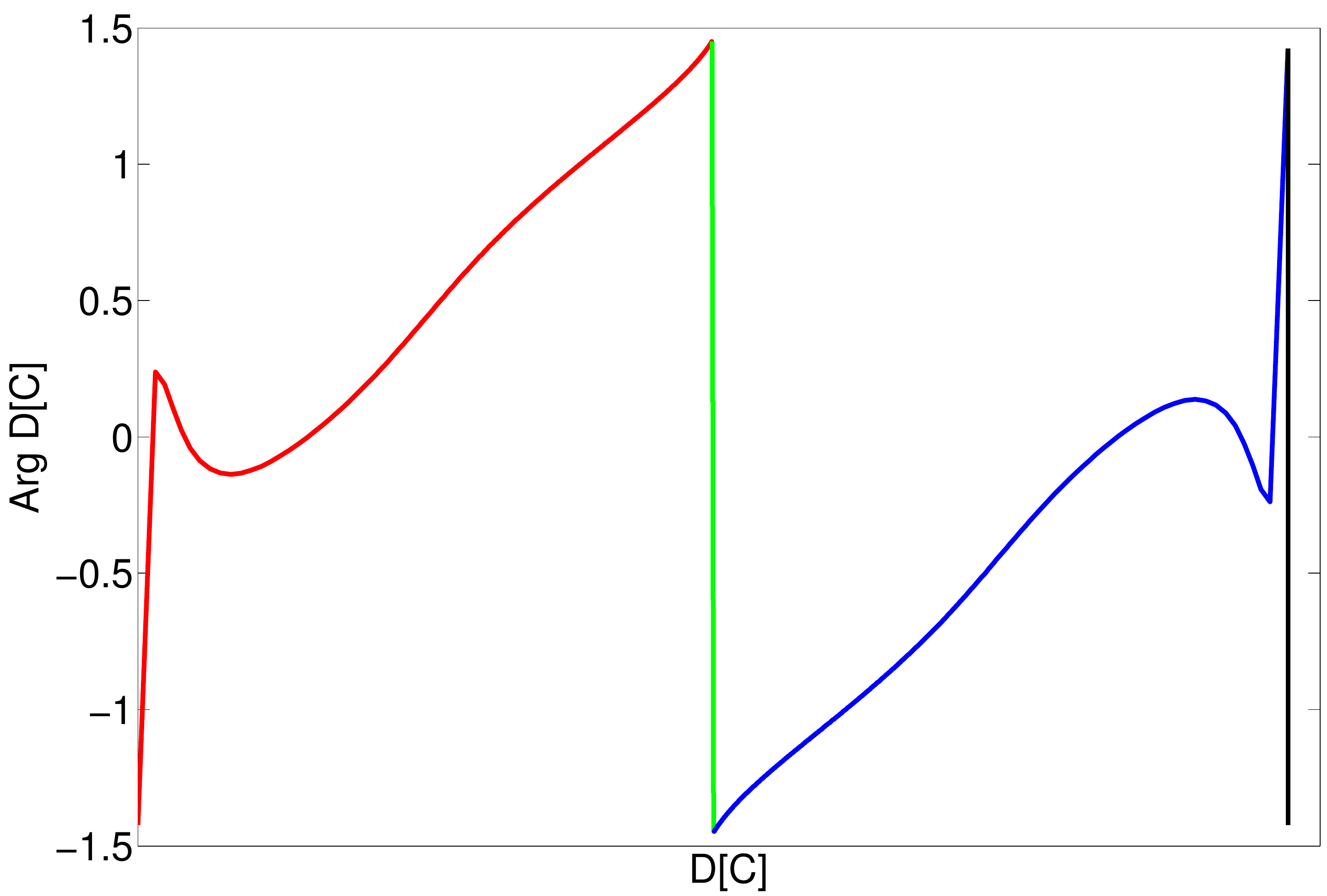}}}
  \caption{Plots $(a)$ and $(b)$ show the change in argument for $D[C]$ corresponding to Figures \ref{evanscontoursmall} and \ref{evanscontourbig} respectively.}
  \label{argchange}
\end{figure}

\noindent We compute the image of $C$ under $D(\lambda)$, which we denote as $D[C]$, and show $D[C]$ in Figure \ref{evanscontoursmall}-\ref{evanscontourbig} for $(r_s,r_b) = (0.1,10)$ and $(0.001,500)$. By the Argument Principle, the number of times $D[C]$ winds around the origin is equal to the number of zeroes of $D(\lambda)$ in the interior of the region enclosed by $C$. Figures \ref{evanscontoursmall}-\ref{argchange} show that the winding number of $D[C]$ around the origin is zero and thus there are no zeroes of the Evans function in the right-half of the complex plane.
\\
\\

\begin{figure}[H]
  \hspace{1mm}\centerline{
  \subfloat[]{\label{evanssmall}\includegraphics[scale=0.3]{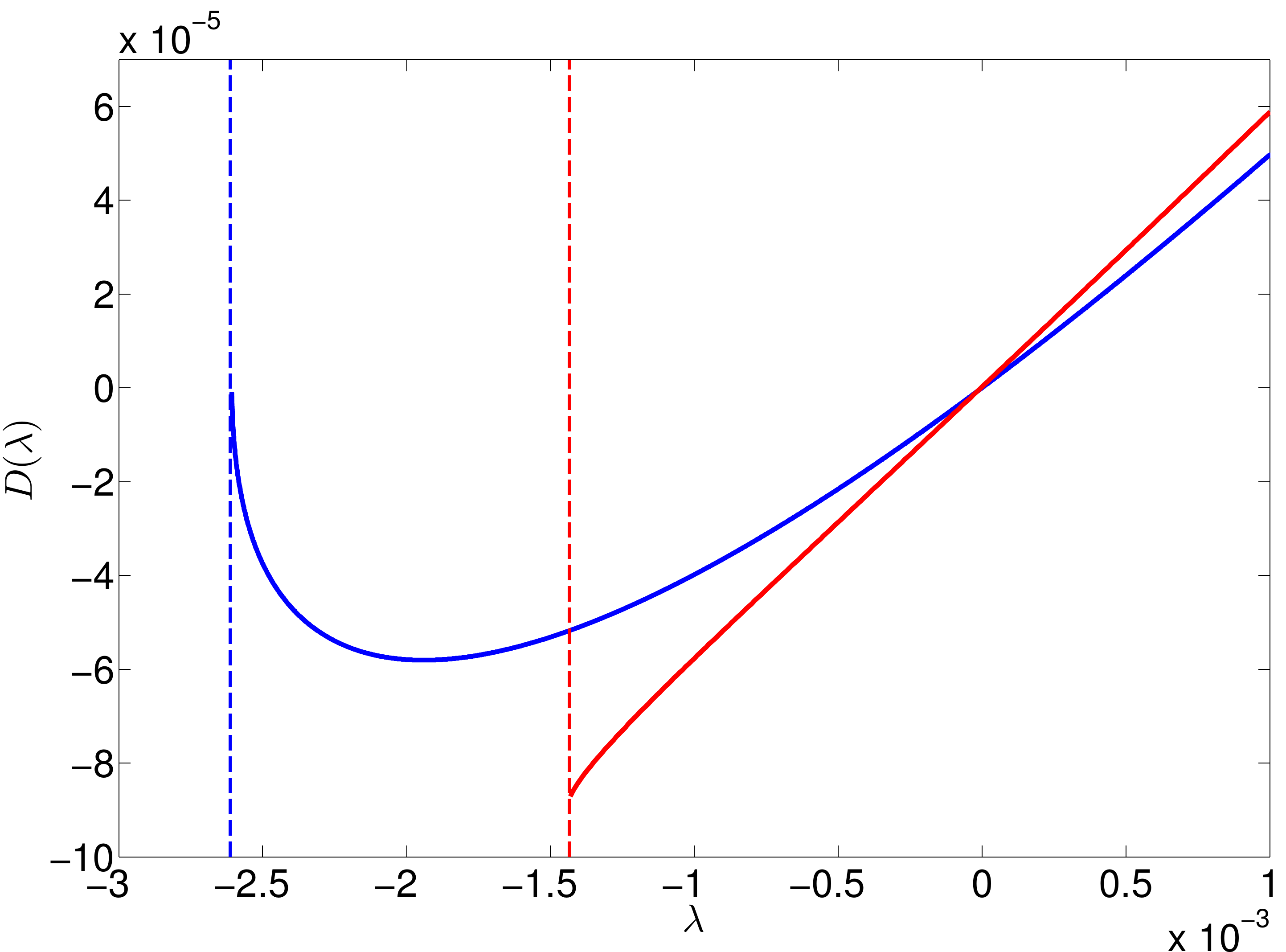}}
  \hspace{1mm}
  \subfloat[]{\label{evansbig}\includegraphics[scale=0.3]{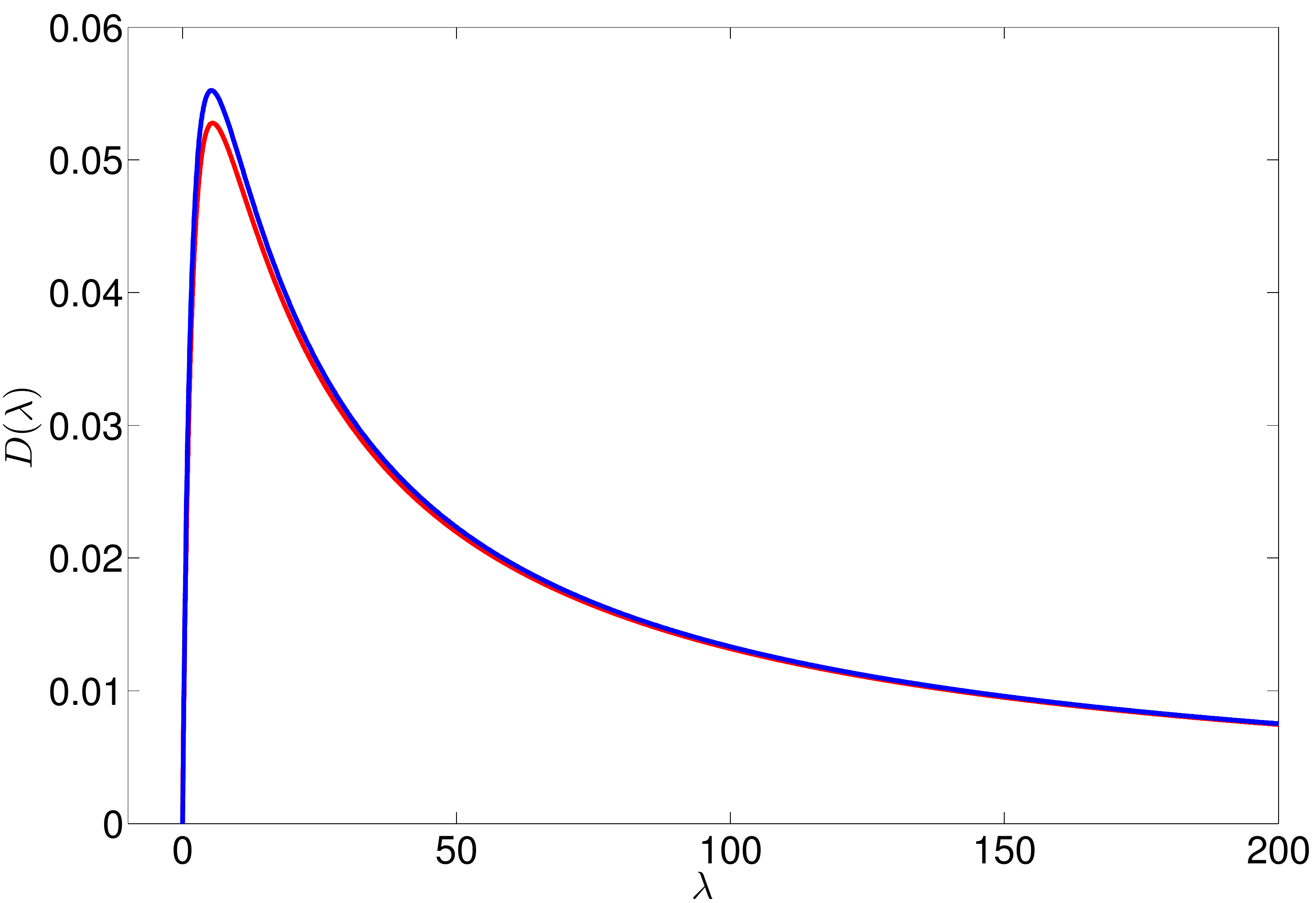}}}
  \caption{Plot of the Evans function given by Eq.~(\ref{evansfunccompound}). The blue and red solid lines in both plots show $D(\lambda)$ with $\alpha=1.1$ and $\alpha=1$ respectively. The dashed lines in plot $(a)$ mark the edge of the absolute spectrum corresponding to each $\alpha$. For $\alpha=1.1$, the only roots are at $\lambda = 0$ and at the edge of the absolute spectrum $\lambda = -0.002607$, whereas for $\alpha=1$ we were unable to detect a zero at the edge of the absolute spectrum due to its proximity to the origin.}
  \label{evansplot}
\end{figure}

\noindent Figure \ref{evansbig} and \ref{evanssmall} show the Evans function computed on the real line for $\lambda \in [0,200]$ and $\lambda \in [-0.002607,0.001]$ respectively. The only roots of the Evans function are at $\lambda = 0$ and $\lambda \approx -0.0026075$, the latter being the edge of the absolute spectrum, which we denote by $\gamma_{A}$. The former is due to translational invariance of the travelling wave solution $(\hat{u},\hat{v})$, while the latter is due to $\gamma_{A}$ being a branch point of $D(\lambda)$. {We note that the root at $\lambda = 0$ is simple, and because of translational invariance must persist throughout all nearby parameter regimes. We have already shown that the root of the Evans function at the branch point is in the left half plane and thus we conclude that no new eigenvalues can be introduced by perturbation.
}
\\
\\
Since system (\ref{PDEsystem}) has no spectrum in the right half plane, the solution given by $\hat{u}$ and $\hat{v}$ is spectrally stable. Moreover, as  the linearised operator $\mathcal{L}$ is an exponentially asymptotic operator, we have that it is also a sectorial operator \cite[see Chapter XVII, $\S$6, Proposition 3 in the former and Example 4.1.8 in the latter]{170,184}. Thus, spectral stability of the travelling wave solution $\hat{u}$ and $\hat{v}$ also implies linear stability. We refer the reader to Section \ref{waveprofilesect} for details on computing $\hat{u}$ and $\hat{v}$.

\noindent We note that for parameter values such that $\gamma_A$ lies closer to the origin, for example $\alpha=1$ (corresponding to $\gamma_A \approx -0.001435$), the method described above fails to detect a zero for the Evans function evaluated at $\gamma_A$, although the qualitative behaviour in the right-half plane remains the same. We show this in Figure \ref{evansplot}.

\section{Wave profile}
\label{waveprofilesect}

\begin{figure}[H]
  \hspace{1mm}\centerline{
  \subfloat[]{\label{waveprofile}\includegraphics[scale=0.3]{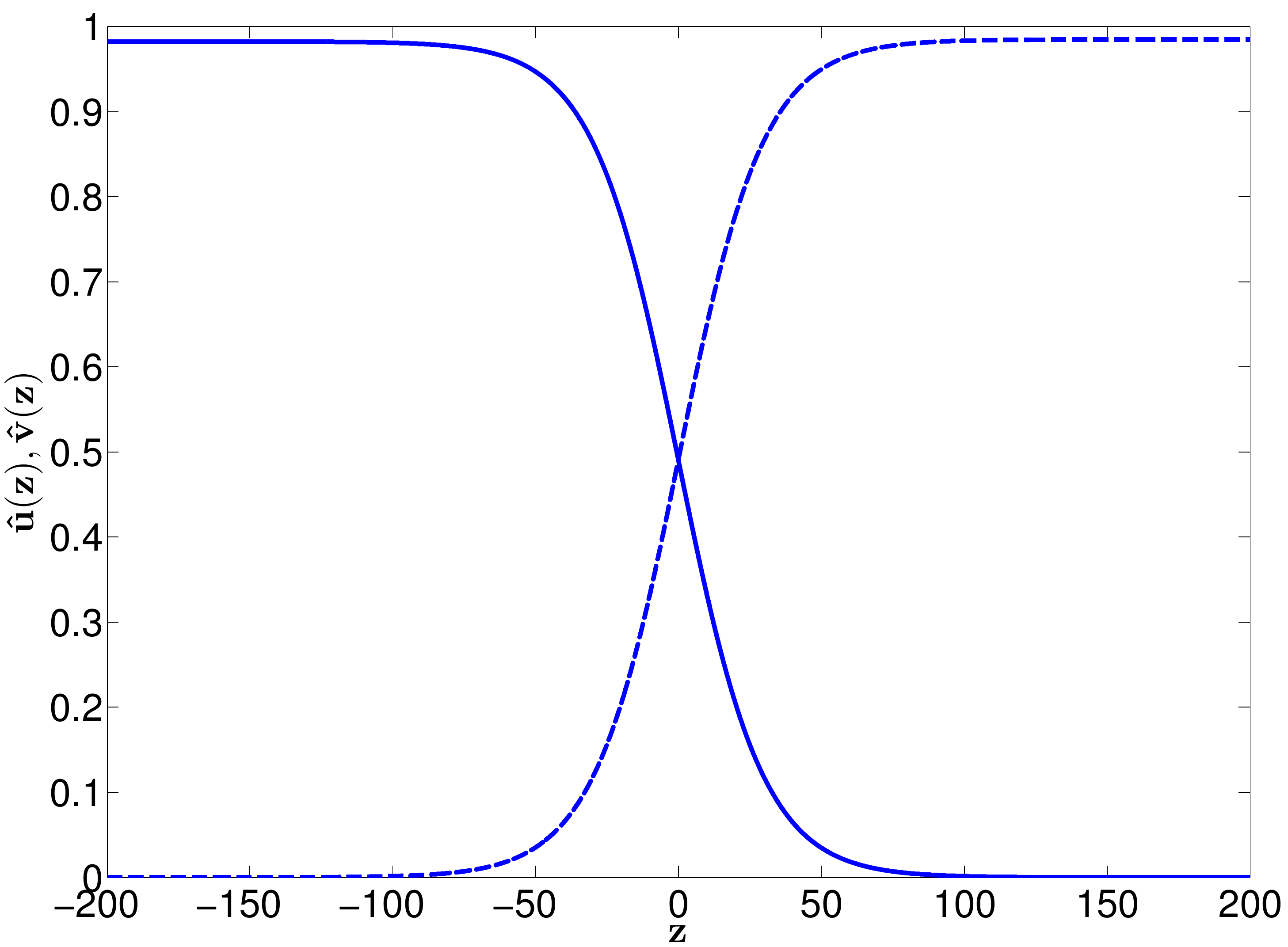}}
  \hspace{1mm}
  \subfloat[]{\label{wavephase}\includegraphics[scale=0.3]{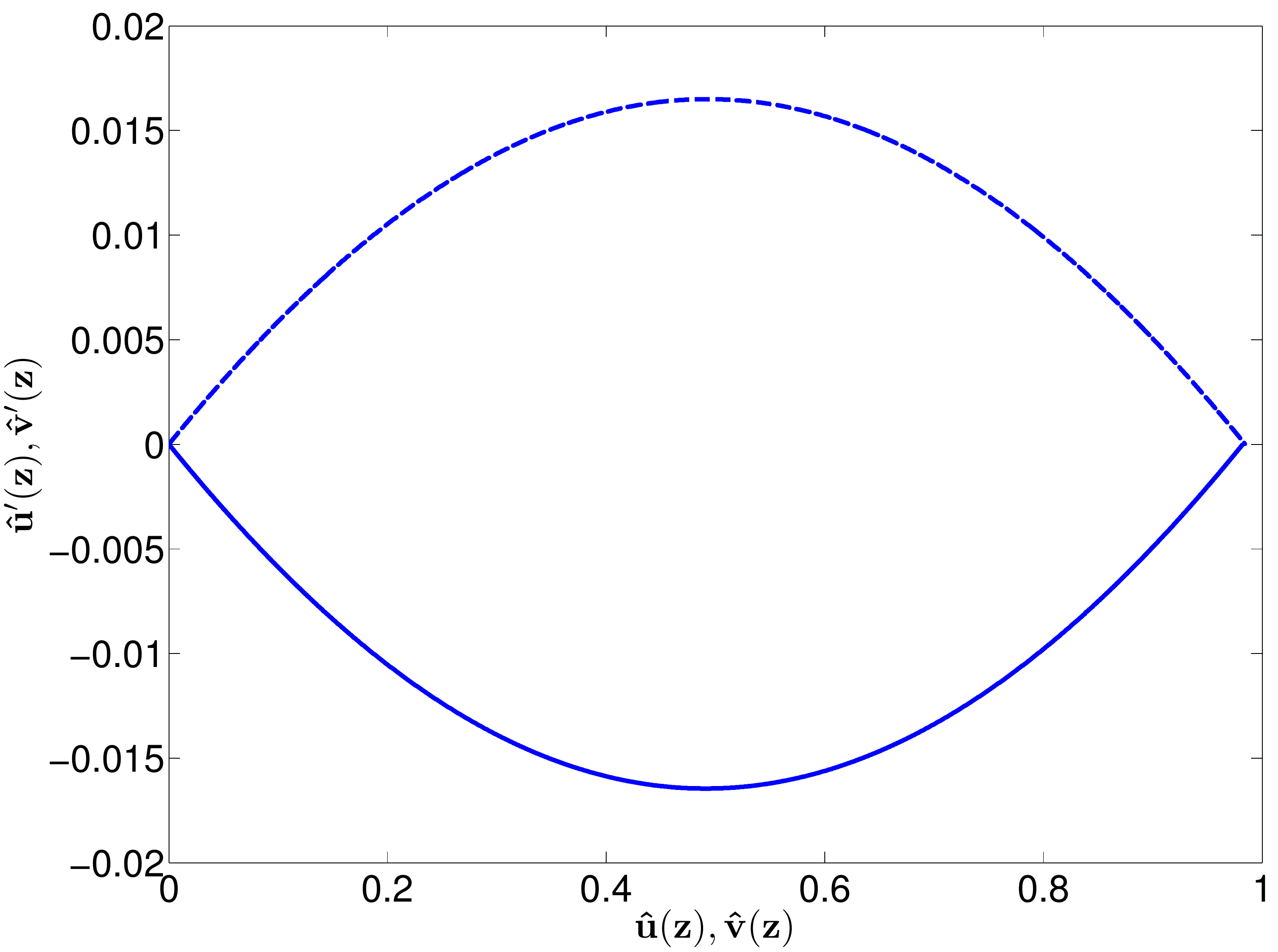}}}
  \caption{Solution to the boundary value problem (\ref{bvp}). Figure $(a)$ shows the wave profile of $\hat{u}(z)$ and $\hat{v}(z)$, represented by solid and dashed lines respectively. Figure $(b)$ shows the heteroclinic connection between equilibrium states $\mathbf{e_{-}}=(1-\alpha \mu,0)$ and $\mathbf{e_{+}}=(0,1-\frac{\mu}{F})$, where the solid and dashed line represent the solution in $u-u'$ and $v-v'$ space respectively.}
  \label{profileplots}
\end{figure}

To compute $A(z,\lambda)$ explicitly at any $z$ requires either a numerical or exact solution for $\hat{u}(z)$ and $\hat{v}(z)$ satisfying Eq.~(\ref{PDEsystem}). We use \textsf{MATLAB}'s \textsf{bvp4c} solver to find a numerical solution corresponding to the case where the populations represented by $\hat{u}$ and $\hat{v}$ are in competition. The boundary conditions listed in (\ref{bvp}) are not sufficient for \textsf{bvp4c} to find a unique solution. We note that since the linearisation of (\ref{travelsys}) as $z \to \pm \infty$ is given by $\mathbf{y}' = A_\infty(0) \mathbf{y}$, we have

\begin{equation}
\label{bc1}
\begin{cases}
\begin{aligned}
u' &\sim \eta^-_2 (u-(1- \alpha \mu))\\
v'  &\sim \eta^-_2 v
\end{aligned}
\end{cases} \quad\quad \text{as } z \to -\infty,
\end{equation}

\begin{equation}
\label{bc2}
\begin{cases}
\begin{aligned}
v' &\sim \eta^+_2 \left( v- \left( 1-\frac{\mu}{F} \right) \right)\\
u'  &\sim \eta^+_2 u
\end{aligned}
\end{cases} \quad\quad \text{as } z \to \infty,
\end{equation}
\\
with $\eta^\pm_2$ as defined in Eq.~(\ref{unstableeigUm})-(\ref{unstableeigUp}), but with $\lambda = 0$. To ensure uniqueness of the solution, we include this information on the derivatives in the boundary conditions by setting

\begin{equation}
\begin{cases}
\begin{aligned}
\frac{u'(-L)}{u(-L)-(1- \alpha \mu)} &= \eta^-_2, \\
\frac{v'(L)}{\left( v(L) - \left( 1-\frac{\mu}{F} \right) \right)} &= \eta^+_2, \\
u  &= \exp(\eta^+_2 L),\\
v  &= \exp(-\eta^-_2 L),
\end{aligned}
\end{cases}
\end{equation}
\\
where $L$ is a large number chosen to represent numerical infinity (chosen to be $L=200$ in Figure \ref{profileplots}). The numerical solution of the wave profiles $\hat{u}$ and $\hat{v}$ are shown in Figure \ref{profileplots}.
\\

\begin{figure}[H]
  \centerline{\includegraphics[scale=0.4]{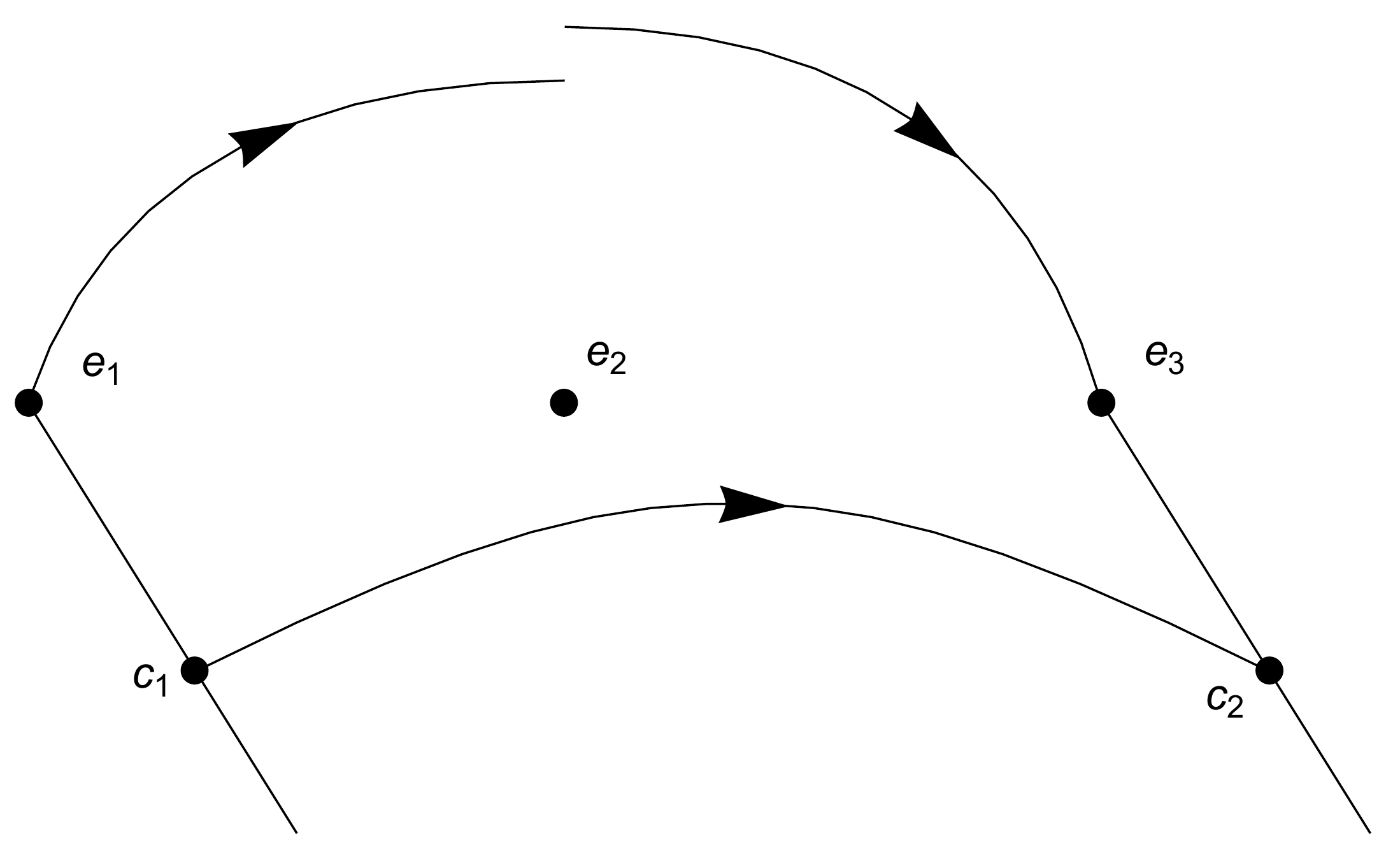}}
  \caption[A diagram illustrating the uniqueness of $c$]{A diagram showing the uniqueness of $c$ by dimension counting, where $e_1$, $e_2$ and $e_3$ denote the three equilibria associated with system \ref{travelsys}. Due to the boundary conditions associated with system \ref{travelsys}, we are only interested in the case where there exists a heteroclinic connection between $e_1$ and $e_3$.}
  \label{prooffigure}
\end{figure}

\noindent Due to the bistability in the system, and from the spatial dynamics of the PDE $u_t = u_{xx} + u(1-u)(u-a)$ (essentially a reaction-diffusion equation with bistable reaction term/strong Allee growth dynamics, see Lewis \& Kareiva \cite{183}), we expect that there is a unique wavespeed $c^*$, which we numerically determine to be approximately 0.027, for which there is a heteroclinic connection between $\mathbf{e}_-$ and $\mathbf{e}_+$ in system (\ref{travelsys}).
\\
\\
Linearising about $\mathbf{e}_-$ and $\mathbf{e}_+$ shows that the dimensions of $W^u(\mathbf{e}_-)$ (the unstable manifold of $\mathbf{e}_-$) and $W^s(\mathbf{e}_+)$ (the stable manifold of $\mathbf{e}_+$) are both equal to 2. We assume that these two manifolds are generic and intersect tranversely, which we interpret as a biologically realistic assumption. Then we consider extending system (\ref{travelsys}) with $c' = 0$. Due to transversality, we have that $ \text{codim} (W^u(\mathbf{e}_-) \cap W^s(\mathbf{e}_+)) = \text{codim}(W^u(\mathbf{e}_-)) + \text{codim}(W^s(\mathbf{e}_+)) = 2+2$, which leads to $\dim(W^u(\mathbf{e}_-) \cap W^s(\mathbf{e}_+)) = 1$. Figure \ref{prooffigure} provides a schematic of this argument.

\section{Discussion}

In this study we have shown the linear stability of a travelling wave solution to a model for \emph{Wolbachia} spread. This is achieved by computing the essential, point and absolute spectrum of the linearised operator and showing the absence of spectrum in the right-half plane. We prove that the essential and absolute spectrum is bounded to the left-half plane for all biologically relevant parameter settings. Due to the numerical nature of locating the point spectrum, we only show that there is no point spectrum in the right-half plane for fixed parameter values. Our results suggest that although \emph{Wolbachia} may be difficult to establish in a local area due to a CI-induced strong Allee effect in the growth dynamics, once it is established the spread of infection is a stable phenomenon.
\\
\\
In addition to our study being an investigation of \emph{Wolbachia} spread dynamics, we present our study as an example of a dynamical systems approach to determining stability of travelling wave solutions in a system of PDEs. Models demonstrating such solutions are ubiquitous, particularly as mathematical modelling is becoming increasingly integrated with the scientific method. The dynamical systems tools we have used in this study can be applied to a wide variety of models currently used in mathematical biology; we believe that their application will improve understanding of the dynamics generated by a model and motivate research in more complicated biological models.
\\
\\
One of the key obstacles impeding the wider use of the tools in this study is the difficulty in computing the point spectrum via the Evans function. There are two key difficulties regarding this. Firstly, numerical methods for evaluating the Evans function sometimes fail, due to the evaluation requiring the solution to a stiff problem. Although in this study we have successfully used the compound matrix method, it is not guaranteed to work for all cases. Secondly, evaluating the Evans function requires the solution whose stability we are interested in. While obtaining such a solution is straightforward when the system of PDEs is exactly solvable, in many cases it is not exactly solvable and instead one must rely on a numerical solution obtained through solving a boundary value problem (see Section \ref{waveprofilesect}). Depending on the dimensionality of the problem and the dimensions of the stable and unstable subspaces at the equilibria, this can be a non-trivial numerical problem.

\bibliographystyle{apalike}
\bibliography{references}
\end{document}